\documentclass[12pt]{article}
\pdfoutput=1 
\usepackage{amsmath}
\usepackage{amssymb}
\usepackage{amstext}
\usepackage{color}
\usepackage{arydshln}
\usepackage{multirow} 
\usepackage{subfig}
\usepackage{graphicx}

\newtheorem{Lemma}{Lemma}

\newtheorem{definition}{Definition}

\newtheorem{theorem}{Theorem}
\newtheorem{corollary}{Corollary}
\newcommand*{\QEDA}{\hfill\ensuremath{\blacksquare}}%
\newtheorem{remark}{Remark}

\newtheorem{assumption}{Assumption}
\newcommand{\rr}{\mathbb{R}}

\newcommand{\dd}{\mathrm{d}}
\newcommand{\scal}[2]{\left\langle #1,#2 \right\rangle}
\newcommand{\norm}[1]{\left\|#1\right\|}
\newcommand{\eig}{\operatorname{Eig}}
\newcommand{\re}{\operatorname{Re}}

\newcommand{\proof}{\textbf{\textit{Proof:}} }
\newcommand{\textwlog}{w\@.l\@.o\@.g\@. }

\begin{document}

\title{Nonlinear stability of flock solutions in second-order swarming models}
\author{J. A. Carrillo, Y. Huang, S. Martin\thanks{ Department of Mathematics, Imperial College
    London, London, SW7 2AZ, UK
    email\textup{\nocorr: \{carrillo,
    yanghong.huang, stephan.martin\}@imperial.ac.uk }
}}

\maketitle

\begin{abstract}
In this paper we consider interacting particle systems which are frequently used to model collective behaviour in animal swarms and other applications. We study the stability of orientationally aligned formations called flock solutions, one of the typical patterns emerging from such dynamics. We provide an analysis showing that the nonlinear stability of flocks in second-order models entirely depends on the linear stability of the first-order aggregation equation. Flocks are shown to be nonlinearly stable as a \emph{family} of states under reasonable assumptions on the interaction potential. Furthermore, we numerically verify that commonly used potentials satisfy these hypotheses and investigate the nonlinear stability of flocks by an extensive case-study of uniform perturbations.  
\end{abstract}



\section{Introduction}

Self-organization and pattern formation are ubiquitous in nature and science, ranging from animal aggregation~\cite{parrish1997animal,couzinKJRF,Couzin2003,cucker2007emergent,cristiani2010modeling} and biological  systems~\cite{Camazine_etal,M&K} to self-assembly of nano-particles~\cite{nano1,nano2}. The intense research during the last two decades are on both individual-based systems and continuum equations~\cite{LR,DCBC,CDMBC,carrillo2010particle}. One of the essential features in these models is the non-zero characteristic speed of the individual agents, which has been modelled in different ways. The speed can be assumed to be constant with a direction based on the averages of the neighbours~\cite{vicsek1995novel} or to be driven by random noise~\cite{degond2008continuum}. On the other hand, a large class of models consist of self-propelled particles powered by biological or chemical mechanisms with friction forces, resulting in a preferred characteristic speed.

In general, the particles with non-zero equilibrium speed do not form any  recognizable patterns~\cite{cucker2007emergent,cfrt}, and interactions within the group have to be included to generate interesting spatial configurations. Most of these interaction forces have been taken into account in the combination of three effects: alignment, repulsion, and attraction; also called the ``first principles of swarming". The basic mechanism account for collisional avoidance and comfort regions (repulsion), grouping and socialization (attraction), and mimetic synchronization (alignment).
The combination of these three effects goes back to fish school modelling \cite{Aoki,Huth:Wissel} and these basic ideas have been improved and applied to several animal species including different mechanisms and interactions more adapted to particular living organisms, see for instance \cite{Ballerini_etal,Cavagna_etal,Hemelrijk:Hildenbrandt,Barbaro_etal,HCH}. 

In this paper,  we focus on the cases with velocity-independent 
 interactions. More precisely, 
by introducing  a pairwise symmetric interaction potential $W(x)=U(|x|)$,
we consider the two-dimensional model
\begin{flalign}\label{micro}
\begin{split}
   \frac{\dd x_i}{\dd t} &= v_i \,, \\
   \frac{\dd v_i}{\dd t} &= \alpha v_i  - \beta v_i \vert v_i\vert^2
   - \sum_{j \neq i} \nabla W(x_i-x_j)\,,
\end{split}
\end{flalign}
where $x_i,v_i\in\rr^2, i=1,\dots,N$ are the positions and velocities of the 
individual particles and $\alpha, \beta$ are
effective values for self-propulsion and friction forces, see
\cite{LR,DCBC,CDMBC} for more discussions.

For this relatively simple system,  a variety of patterns are observed,
for instance coherent moving flocks and rotating mills~\cite{LR}. 
Other patterns like rings and clumps were discovered in~\cite{DCBC,CDMBC},
with the introduction of the concept of \emph{H-stability} from 
statistical physics to better characterize the patterns in the parameter space.
Delta rings, uniform distribution on a circle, have been further studied thoroughly for power-law like 
potentials, we refer to \cite{KSUB,predict,soccerball,BCLR,BCLR2,bigring,ABCV} 
for details. While the stability and bifurcation of the ring solutions
can be investigated in a straightforward way because of their explicit 
uniform particle representations on a circle, 
there are few studies on the more prevalent 
compact steady solutions, like flocks or mills. The reason lies
in the difficulties to solve some complicated integro-differential equations
for most of the popular choices of the potential $W$. Only for certain particular
potentials like the quasi-Morse type proposed in~\cite{Carrillo2013},  
the existence and uniqueness of the coherent moving flock can be
established rigorously \cite{CHM}.

In this paper, we focus on the stability of flock solutions 
of~\eqref{micro} for general potentials defined below. 
\begin{definition}
A flock solution of the particle model \eqref{micro} is a spatial configuration $\hat{x}$ with zero net
interaction force on every particle, that translates at a uniform
velocity $m_0\in\rr^2$ with $|m_0|=\sqrt{\frac{\alpha}{\beta}}$, hence $(x_i(t),v_i(t))=(\hat{x_i}-tm_0,m_0)$. 
\label{def-flock} 
\end{definition}

We note, that the spatial configuration $\hat{x}$ is a
stationary state to the first-order interacting particle system
\begin{equation} \label{aggregation}
\frac{\dd x_i}{\dd t} = - \sum_{j \neq i} \nabla W(x_i-x_j)\,,  
\end{equation} 
which is analysed for example in \cite{nano1,nano2}. In this work, we
focus on general flock solutions whose spatial configurations
tend to a compactly supported particle density in the continuum
limit.

One conclusion from our analysis is the somehow
unexpected deep relation between the linear stability of the flock
spatial configuration as steady state for the first-order swarming
model \eqref{aggregation} and the nonlinear stability of the family of
associated flock solutions for the second-order swarming model
\eqref{micro}. This relation was already found in the linear
stability analysis around flock solutions in \cite{ABCV}. There the authors showed that the linearization of
\eqref{aggregation} around the equilibrium state $\hat{x}$ has a positive eigenvalue if and only if
the linearization of \eqref{micro} around the steady flock solution in the co-moving frame has an eigenvalue with positive real part. Moreover, they show that if the equilibrium state $\hat{x}$ is linearly stable for
\eqref{aggregation}, then the associated flock solution is always linearly unstable due to the presence of a generalized eigenvector associated to the zero eigenvalue of the linearization of the flock solution in the co-moving frame due to symmetries.

However, it is in this work where we clarify completely this deep
imbrication by showing the nonlinear stability of the flock solutions under mild hypotheses on the linearization of \eqref{aggregation} about the spatial configuration $\hat{x}$. 
Actually, the basic hypotheses, apart from few technical assumptions, are that
the steady state $\hat{x}$ of \eqref{aggregation} is linearly asymptotically
stable except the obvious symmetries: translations and rotations. Therefore, we can apply
our theorem to verify the nonlinear stability of flock solutions to more
biologically adapted potentials such as the Morse and the Quasi-Morse potentials
\cite{DCBC}.

The main result of this paper asserts, except few technicalities,
that \emph{the family of flock solutions} associated to a linearly
asymptotically stable steady state of \eqref{aggregation} is
\emph{asymptotically stable} for the dynamics of \eqref{micro}. 
Here, the asymptotic stability
of the family of flock solutions means that any small enough
perturbation in $(x,v)$-space at any time $t_0$ will, 
under the dynamics of the system \eqref{micro},
relax towards (likely) another flock solution in the family at an
exponential rate as $t\rightarrow\infty$. Let us finally emphasize
that the most rigorous way of stating our main theorem uses
advanced concepts of dynamical systems \cite{HPS}. Our main theorem can be rephrased
as follows: the family of flock solutions to \eqref{micro}
associated to a linearly (except symmetries) asymptotically stable steady state of
\eqref{aggregation} forms a normally hyperbolic invariant
manifold for the system \eqref{micro} with an empty unstable manifold.

This result has important applications, especially in the study 
of flock patterns of~\eqref{micro} using particle simulations.
In general, the desired patterns are observed only after certain
transit dynamics, and with carefully prepared initial data. With the 
spatial configurations from the first order system~\eqref{aggregation} 
(if they exist), the computational load is reduced significantly,
making possible to study the existence of these patterns in great
detail.

After some preliminary notations and assumptions, the main result of this paper is stated in Section 2. The formulation of the linearized dynamics when the mean velocity is transformed to zero is given in Section 3, where the non-trivial relation between the Jacobian matrices of the two systems \eqref{micro} and \eqref{aggregation} is clarified. In Section 4, we use the eigen-space structures of these Jacobian matrices to reveal the connection between the nonlinear stability of the flock solutions for \eqref{micro} and the linear stability of their spatial configuration for \eqref{aggregation}, leading to the proof of the main theorem. Section 5 is devoted to numerical experiments to check the validity of the assumptions on the linear stability for the steady state of the system \eqref{aggregation} for different biologically reasonable potentials.

%%%%%%%%%%%%%%%%%%%%%%%%%%%%%%%%%%%%%%%%%%%%%%%%%%%%%%%%%%%%%%%%%%%%%%

\section{Main Result}
In this section, we formulate the main theorem of this work.
We will start with some assumptions on the spatial
configuration of flock solutions, which is a stationary state of the first-order swarming
model \eqref{aggregation}. The hypotheses assume this spatial configuration
to be stable under the dynamics of \eqref{aggregation}. 

Given a steady spatial particle configuration to
\eqref{aggregation}, i.e, a set of particle positions
$(\hat{x}_i)_{i=1}^N$ such that
\begin{equation}\label{steadyparticle}
\sum_{j \neq i} \nabla W(\hat{x}_i-\hat{x}_j)=0\,,\qquad \mbox{
for all } i\in\{1,\dots,N\}\,,
\end{equation}
we denote by $\hat{x}$ the $2N$-dimensional spatial configuration
vector $(\hat{x}_1,\dots,\hat{x}_N)$. We can now write the
linearization of the system \eqref{aggregation} at
the steady configuration $\hat{x}$ as
\begin{equation}\label{linearizedaggr}
\frac{\dd h}{\dd t} =  G(\hat{x})h\,,
\end{equation}
where $G(\hat{x})$ denotes the Jacobian matrix associated to
\eqref{aggregation} and $h\in \rr^{2N}$. The $2N\times
2N$-Jacobian matrix is explicitly given by $G(\hat{x}) = (G_{ij})$
with $G_{ij}$ being the $2\times 2$-blocks defined as
\begin{equation}\label{hessian}
G_{ij}=\begin{cases} -\displaystyle \sum_{j \neq i} \mbox{Hess }
W(\hat{x}_i-\hat{x}_j) & \mbox{for } i=j\\[2mm] \,\,\,\,\,\,\,  \mbox{Hess } W(\hat{x}_i-\hat{x}_j) & \mbox{for } i\neq j
\end{cases}\,,
\end{equation}
with $\mbox{Hess } W$ denoting the Hessian matrix of $W$.

We will make extensive use of the following standard notation:
given two matrices $A\in\rr^{n\times m}$ and
$B\in\rr^{p\times q}$, the Kronecker product $A\otimes B$ is
defined as the matrix $A\otimes B = (a_{ij}B)_{ij} \in
\rr^{np\times mq}$. Let us also denote by $0_n, 1_n$ the column
vectors of length $n$ with all entries equal to $0$ or $1$ respectively, 
and $I_n$ the identity matrix of size $n$.
Let us now specify the precise assumptions on $\hat{x}$, and thus
indirectly on the potential $W$.

\begin{assumption}
Given an interaction potential $W$, we assume the existence of a
linearly asymptotically stable (except symmetries) stationary state not lying on a
line, i.e., a spatial configuration $\hat{x}$ such that the
following assumptions hold:
\begin{itemize}
\item {\rm\bf (H1)} $\hat{x}$ is a stationary state of
$\eqref{aggregation}$, i.e. $\hat{x}$ satisfies \eqref{steadyparticle}.

\item {\rm\bf (H2)} The eigenspace for the zero eigenvalue of
$G(\hat{x})$ is spanned by the eigenvectors
\begin{align}\label{evdef-w3}
w_1\!=\! 1_N\!\otimes\!\begin{pmatrix}
1 \\ 0 \end{pmatrix},\ 
w_2\!=\!1_N\!\otimes\!\begin{pmatrix} 0\\ 1\end{pmatrix}\!, \mbox{ and }
w_3\!=\!\left(\!
 I_{N} \otimes
 \begin{pmatrix}
 0 &-1 \\ 1 &0
\end{pmatrix}\!\right)\! \hat{x}
\end{align}
representing invariance of $G(\hat{x})$ with respect to
translations and rotations in $\rr^2$, and therefore is three dimensional, $\dim(\eig(G(\hat{x}),0)) = 3$.

\item {\rm\bf (H3)} All other eigenvalues of $G(\hat{x})$ have
negative real parts, i.e. $\re(\lambda_i)<0 \Leftrightarrow
\lambda_i\neq 0$ for all $i\in\sigma(G(\hat{x}))$.

\item {\rm\bf (H4)} Not all points of $\hat{x}$ lie on a straight
line in $\rr^2$.
\end{itemize}
\end{assumption}

Hypotheses \textbf{(H1)--(H3)} imply that $\hat{x}$ is linearly
asymptotically stable for \eqref{linearizedaggr} except for the
obvious symmetries, or equivalently that the Jacobian matrix
$G(\hat{x})$ is negative definite except for the null eigenspace
in \eqref{evdef-w3}. As a consequence, $\hat{x}$ is nonlinearly
asymptotically stable for the dynamics of \eqref{aggregation}
except for translations and rotations. More rigorously speaking 
and making use of more advanced dynamical systems theory \cite{HPS}, 
properties \textbf{(H1)--(H3)} imply that the family of stationary 
states of \eqref{aggregation} obtained by rotations
and translations from $\hat{x}$ given by
\[
RT(\hat{x})=\{y: y=1_N\otimes b+(I_N \otimes R(\phi))\hat{x}\,,\, 
\phi\in[0,2\pi), b\in\rr^2\}
\]
forms a normally hyperbolic invariant manifold for 
the system of ODEs in \eqref{aggregation} for which the unstable manifold 
is empty. Here, $$R(\phi)=\begin{pmatrix} \cos\phi & -\sin \phi \\
\sin\phi &\cos\phi \end{pmatrix}$$ is the $2\times 2$ rotation matrix of 
angle $\phi$.  This stability concept will be further discussed 
and used later on for the second order model \eqref{micro}.

Let us remark, that assumptions \textbf{(H1)--(H4)} are very mild
in the sense that they are natural properties of stationary states of $\eqref{aggregation}$ found through
particle simulations from generic initial configurations, e.g.
compactly supported flocks. Particular special cases satisfying
hypotheses \textbf{(H1)--(H4)} are the so-called Delta rings:
particles arranged equidistantly on a circle. 
We will discuss in the last section
more biologically relevant potentials, such as Morse-like potentials
\cite{nano1,nano2,LR,DCBC,Carrillo2013,CHM}, for which we can numerically find stationary states
$\hat{x}$ satisfying hypotheses \textbf{(H1)--(H4)}.

Now, let us make precise the concept of the family of flock solutions
of \eqref{micro} associated to the stationary state $\hat{x}$ of
\eqref{aggregation}. Let us denote by $\mathcal{V}_s$ the set of possible
steady velocity vectors for a flock solution of speed
$s=\sqrt{\tfrac{\alpha}{\beta}}$ to \eqref{micro} defined as
$\mathcal{V}_s=\{m_0\in\rr^2: |m_0|=s\}$. Then the family of flock
solutions of \eqref{micro} associated to the steady state
$\hat{x}$ of \eqref{aggregation} is given by
\begin{equation}
\begin{gathered}
\mathcal{Z}_F = \left\{z= \begin{pmatrix}
x^*+v^* t \\
v^*
\end{pmatrix}_{t\in\rr}, v^*=1_N\otimes m_0,
x^*\in  RT(\hat{x}), m_0\in \mathcal{V}_s \right\}\,.
\end{gathered}\label{flocks}
\end{equation}

It is straightforward to check that all $z\in\mathcal{Z}_F$ are
solutions to \eqref{micro} in vectorial form. The next section will be devoted to
analyse the system \eqref{micro} in suitable coordinates to study their stability. 
A change of variables to \eqref{micro} that eliminates the translation of one particular
member of the flock family is given by  the co-moving frame $\tilde{v}_i := v_i - m_0$
transforming the flock $(x^*+v^* t, v^*)$ onto a stationary
state. However, it is very natural to think of small perturbations
on $(x^*+tv^*, v^*)$ that, under the dynamics of \eqref{micro},
tend to possibly another spatial configuration in $RT(\hat{x})$
translating at a different velocity generated by $\tilde{m}_0$,
i.e. to another member of the family of flocks $\mathcal{Z}_F$.
Hence, we cannot expect linear stability for general perturbations
for such a change of variables, see \cite{ABCV} where the linear
instability of a flock solution in the co-moving frame is shown.
Instead, we propose a new change of variables
measuring the difference to mean velocity in velocity space, and
aim to establish nonlinear stability in the \emph{entire family}
$\mathcal{Z}_F$.

\begin{definition}
Defining the mean velocity of $N$ particles as
\begin{equation*}
m(t) := \frac{1}{N}\sum_{i=1}^{N} v_i(t)\,,
\end{equation*}
we introduce the change of variables with respect to the mean velocity
$m(t)$ as:
\begin{align}
\tilde{x}_i(t) &:= x_i(t) - \int_0^{t} m(s) \dd s\,,\nonumber\\
\label{newcoord}\\[-3mm]
\tilde{v}_i(t) &:= v_i(t) - m(t).\nonumber
\end{align}
\end{definition}

We will work with the dynamical system \eqref{micro} and its
linearization in vector form, and thus we need to introduce some
matrix notation used in the following.
\begin{definition}[Notation]
\begin{enumerate}
\item $0_{n\times p}$ is a block of zeros of size $n\times p$.

\item Writing $\delta^i_n\in\rr^n$, we refer to a vector of length
$n$ with value $1$ as its i-th entry and zeros elsewhere.

\item For the sake of simplicity, we occasionally write column
vectors as row vectors neglecting the transpose operation on its
components, e.\@ g.\@
  \begin{equation*}
 (x,v,0_2)^T:=
  (x^T,v^T,0_2^T)^T=
 \begin{pmatrix}
  x\\ v \\ 0_2
  \end{pmatrix}.
    \end{equation*}

\item For a square matrix $A\in\rr^{n\times p}$ or a column vector
$b\in\rr^{n}$, we denote by
  \begin{equation*}
 \lceil\! A\rceil := (a_{ij})_{\substack{i=1,\dots n-2 \\ \!\!\!\!\!\! 
 j=1,\dots,p}}  \quad,\quad  
 \lceil b\rceil= (b_i)_{i=1,\dots,n-2}.
  \end{equation*}
the matrix obtained from $A$ by eliminating the last two rows and
analogously for vectors. The $k$-th row or column of $A$ is
labelled as $a_{k,\cdot}$ or $a_{\cdot,k}$ respectively.
  \end{enumerate}
  \label{def-notations}
  \end{definition}

\begin{Lemma} Using the change of coordinates \eqref{newcoord}, the microscopic
model \eqref{micro} reads
  \begin{align}
  \frac{\dd \tilde x_i}{\dd t} &= \tilde v_i, \nonumber\\
  \frac{\dd \tilde v_i}{\dd t} &= \frac{N\!-\!1}{N}\left(\alpha-\beta|\tilde v_i+m|^2\right)(\tilde v_i+m) - \frac{1}{N}\sum_{j\neq i} \left(\alpha-\beta|\tilde v_j+m|^2\right)(\tilde v_j+m) \nonumber\\
  &\quad - \sum_{j\neq i}\nabla W(\tilde x_j-\tilde x_i),  \label{newcoord-veq}\\
  \frac{\dd m}{\dd t} &=\frac{1}{N} \sum_{j}
  \left(\alpha-\beta|\tilde v_j+m|^2\right)(\tilde v_j+m).\nonumber
  \end{align}
For any mean velocity $m_0\in \mathcal{V}_s$, the vector
$\hat{Q}=(\hat{x}, 0_{\,2N} , m_0)^T $ is a stationary state of
\eqref{newcoord-veq}. Moreover, the Jacobian matrix $F$ obtained
by linearizing \eqref{newcoord-veq} around a stationary state
$\hat{Q}$ with an arbitrary orientation $m_0$ is given by
\begin{equation}
F\!=\!  \begin{pmatrix}
  0_{2N\times 2N} & I_{2N} & 0_{2N\times 2} \\[2mm]
  G(\hat{x}) & \!\!\!- \begin{pmatrix}
  \frac{N-1}{N} & -\frac{1}{N} & \dots &-\frac{1}{N} \\
  -\frac{1}{N} & \ddots & &\vdots \\
  \vdots & & \ddots & -\frac{1}{N} \\
  -\frac{1}{N} & \dots & -\frac{1}{N} & \frac{N-1}{N}
  \end{pmatrix}\! \otimes \! 2\beta (m_0\otimes m_0^{T})
  &
  0_{2N \times 2} \\[11mm]
  0_{2\times 2N} & - 1^T_{N}\otimes \frac{2\beta}{N}(m_0\otimes m_0^{T}) & \!\!\!\!\!\!\!\!\!\!\!\!\!\!\!\!\!\!\!\!-2\beta (m_0\otimes m_0^{T})
  \end{pmatrix}
  \label{Jacobianofmicronewcoord}
\end{equation}
\proof
By definition $\frac{\dd \tilde{x}_i}{\dd t} = \frac{\dd x_i}{\dd t} - m(t) = v_i(t) - m(t) = \tilde{v}_i(t). $ Next we note $\tilde{x}_i-\tilde{x}_j  = x_i - x_j$, and hence $\nabla W (\tilde{x}_i -\tilde{x}_j ) = \nabla W(x_i - x_j) $. Now, we compute
\begin{align*}
\frac{\dd m}{\dd t} &= \frac{1}{N}\sum_i \frac{\dd v_i}{\dd t}
= \frac{1}{N}\left( \sum_i \left(\alpha-\beta |v_i|^2\right)v_i  + \sum_i \sum_{j\neq i}  
\nabla  W(x_i - x_j) \right)\\
&= \frac{1}{N} \sum_i \left(\alpha-\beta |\tilde{v}_i  + m|^2\right)(\tilde{v}_i +m)\,,
\end{align*}
using that the net total force in the system is zero by symmetry, i.e. Newton's action-reaction principle. From this, the equation for 
$
\tfrac{\dd \tilde{v}_i }{\dd t}=\tfrac{\dd v_i}{\dd t}- \tfrac{\dd m}{\dd t}
$ is readily available.
Now, any state  $\hat{Q}=(\hat{x}, 0_{\,2N} , m_0)^T$, with $m_0\in \mathcal{V}_s$ is a stationary state of \eqref{newcoord-veq} since $G(\hat{x})=0$ and all self-propulsion terms vanish.
Let us drop the tildes henceforth for simplicity and let us compute the Jacobian matrix of the system \eqref{newcoord-veq}. The first line of blocks in $F$ is trivial. In the second line of blocks, we first have $G(\hat{x})$ since the dependency on $x$ in \eqref{newcoord-veq} is identical to \eqref{micro}, hence one obtains the Jacobian of the first-order particle system \eqref{aggregation} evaluated at the stationary spatial configuration.
To compute the second block, consider an arbitrary particle $i$ and its velocity vector $(v_{i,1}, v_{i,2})^T$. Then, identifying the right-hand sides of \eqref{newcoord-veq} with the time derivative symbols, we get:
%$
\begin{gather*}
\partial_{v_{i,1}} \frac{\dd v_{i,1}}{\dd t} = \frac{N-1}{N}\left( -2\beta(v_{i,1}+m_{1})^2 + (\alpha-\beta|v_i+m|^2)   \right),
%$
\\
%$
\partial_{v_{i,2}} \frac{\dd v_{i,1}}{\dd t} = \frac{N-1}{N}\left( -2\beta(v_{i,2}+m_{2})(v_{i,1}+m_1) + (\alpha-\beta|v_i+m|^2)   \right),
%$
\\
%$
\partial_{v_{j,1}} \frac{\dd v_{i,1}}{\dd t} = -\frac{1}{N} \left( -2\beta (v_{j,1}+m_1)^2 + (\alpha - \beta|v_j+m|^2)\right) ,
%$
\end{gather*}
for $j\neq i$ and analogous expressions for $\partial_{v_{i,1}} \frac{\dd v_{i,2}}{\dd t}, \partial_{v_{i,2}} \frac{\dd v_{i,2}}{\dd t}, \partial_{v_{j,2}} \frac{\dd v_{i,1}}{\dd t},   \partial_{v_{j,1}} \frac{\dd v_{i,2}}{\dd t} $, and $ \partial_{v_{j,2}} \frac{\dd v_{i,2}}{\dd t}$.
Evaluating at $\hat{Q}$ the last expression in all terms vanishes since $|m_0|^2=\tfrac{\alpha}{\beta}$, which in total gives the structure of the block using the Kronecker product notation.
In the third line of $F$, the second block follows analogously. For the last block, we compute
$$
\partial_{m_1}\frac{\dd m_1}{\dd t} = \frac{1}{N} \sum_{j} \left(-2\beta (v_{j,1}+m_{1})^2 +(\alpha-\beta|v_{j}+m|^2)\right)\,,
$$
which reduces to
$
\partial_{m_1}\frac{\dd m_1}{\dd t} = -2\beta m_1^2
$ at $\hat{Q}$, and analogously for the remaining partial derivatives.
This completes the computation of the structure of $F$.
\QEDA
\label{cor-newcoord}
\end{Lemma}

Let us point out that the family of flock solutions to \eqref{micro} is translated through the change of variables \eqref{newcoord} to the set of stationary solutions of \eqref{newcoord-veq} given by
\begin{equation}
\begin{gathered}
\tilde{\mathcal{Z}}_F = \left\{Q^*=\begin{pmatrix}
x^*\\
0_{2N}\\
m_0
\end{pmatrix}, \,x^*\in  RT(\hat{x}), m_0\in \mathcal{V}_s \right\}\,.
\end{gathered}\label{flocks2}
\end{equation}

Now, we have almost all the ingredients to write the main result of this paper that holds under an additional technical assumption on the interplay between the linearization and the orientation of the flock.

\begin{assumption}[H5] All eigenvectors to non-zero eigenvalues of $G(x^*)$ are not particle-wise orthogonal to the fixed mean-velocity $m_0$. This means, that for pairs of entries $2i-1,2i$ of any eigenvector $w$ of $G(x^*)$, we have
\begin{equation*}
\scal{m_0}{w_{\substack{2i\text{-}1\\ \!\!\!\!\!2i}}} \neq 0 \text{ for at least one } i\in{1,\dots,N}.
\end{equation*}
\end{assumption}

We will understand the role of this additional assumption below, but let us point out that this 
assumption is generically met by the members of the flock family $\tilde{\mathcal{Z}}_F$. The main result of this work reads as:

\begin{theorem}\label{theo-main}
Assume that the symmetric pairwise interaction potential $W$ allows for a stationary spatial
configuration $\hat{x}$ of \eqref{aggregation} satisfying assumptions {\rm \bf (H1)-(H4)}.
Consider the second-order swarming model $\eqref{micro}$ with
$\alpha,\beta>0$, then the family of flock solutions $\mathcal{Z}_F$ of \eqref{micro} defined in \eqref{flocks} is {\rm locally asymptotically stable} for all configurations $x^*$
satisfying {\rm \bf (H5)} in the following sense: Any small enough
perturbation in $(x,v)$-space at any time $t_0$ of the flock
solution $z^*\in \mathcal{Z}_F$ associated to $x^*$ will, under the dynamics of the
system \eqref{micro}, relax towards another flock solution
$\tilde{z}\in \mathcal{Z}_F$ at an exponential rate as
$t\rightarrow\infty$. 
\end{theorem}

Let us remind that assumptions {\rm \bf (H1)-(H4)} imply in particular that the family of stationary states $RT(\hat{x})$ forms a normally hyperbolic invariant manifold for \eqref{aggregation} with empty unstable manifold. Therefore, our main theorem can be rephrased as stating that the family of flock solutions associated to $\hat{x}$ of the second order model \eqref{micro} has essentially the same property. More precisely, we have the following result.

\begin{corollary}\label{cor-main}
Assume that the symmetric pairwise interaction potential $W$ allows for a stationary spatial
configuration $\hat{x}$ of \eqref{aggregation} satisfying assumptions {\rm \bf (H1)-(H4)}.
Then the subset of stationary states of $\tilde{\mathcal{Z}}_F$ in \eqref{flocks2} satisfying {\rm \bf (H5)} is a normally hyperbolic invariant manifold of \eqref{newcoord-veq} with empty unstable manifold.
\end{corollary}

Next two sections are devoted to the proof of these results. We start in next section by restricting the dynamics to a physically consistent subspace.

%%%%%%%%%%%%%%%%%%%%%%%%%%%%%%%%%%%%%%%%%%%%

\section{Physically Consistent Dynamics and its Linearization}

Note that \eqref{newcoord-veq} is a $4N+2$ dimensional system and thus $F\in\rr^{(4N+2) \times (4N+2)}$ by having added the mean velocity as an additional variable. However, a stability analysis of \eqref{Jacobianofmicronewcoord} would include unphysical perturbations where particle velocities do not match mean velocity. The dynamics we are interested in are $4N$-dimensional, and thus we need to reduce \eqref{Jacobianofmicronewcoord} to physically consistent dynamics. Clearly, mean velocity consistent states are invariant under the dynamics of \eqref{newcoord-veq}, since
\begin{equation*}
\frac{\dd}{\dd t}\frac{1}{N} \sum_i\left(v_i - m\right) = 0.
\end{equation*}
We hence define a base for the $4N$-dimensional subspace of mean velocity consistent states by removing the redundant equation of the $N$-th particle.

\begin{definition}
Define the subspace of mean velocity consistent states as $\operatorname{span}(B)$ with
$
B = \left\{
b_i \right\}_{i=1,\dots,4N}
$
and
\begin{equation*}
b_i = \begin{cases}
\begin{pmatrix} \delta_{2N}^i,& 0_{2N},& 0_2 \end{pmatrix}^T & i=1,\dots,2N, \\
\begin{pmatrix}0_{2N},& \delta_{2N-2}^k &-1 &0,& 0_2 \end{pmatrix}^T & i=2N+(2k-1)\,,\, k=1,\dots,N-1, \\
\begin{pmatrix} 0_{2N}, &\delta_{2N-2}^k &0& -1,& 0_2\end{pmatrix}^T & i=2N+2k\,,\, k=1,\dots,N-1 ,\\
\begin{pmatrix} (0_{2N}, & 0_{2N},& \delta_2^k\end{pmatrix}^T & i=4N-2+k \,,\, k=1,2.
\end{cases}
%\label{base}
\end{equation*}
Any physically relevant perturbation of \eqref{micro} is projected onto $\operatorname{span}(B)$ via
\begin{equation*}%\label{projection}
P : \begin{pmatrix}
\Delta{x} \\ \Delta{v}
\end{pmatrix} \mapsto \begin{pmatrix}
\Delta{x} \\[1mm] \hdashline \\[-3mm] \Delta{v}_1 - \frac{1}{N}\sum\Delta{v}_i \\[1mm] \vdots \\[1mm] \Delta{v}_N - \frac{1}{N}\sum\Delta{v}_i \\[1mm] \hdashline \\[-3mm]  \frac{1}{N}\sum\Delta{v}_i
\end{pmatrix} \in \operatorname{span}(B).  
\end{equation*}
\end{definition}

Clearly, $B$ is isomorphic to a base of $\rr^{4N}$ and the relevant invariant subspace for the dynamics of \eqref{newcoord-veq}. As mean velocity consistency is a linear algebraic constraint to \eqref{newcoord-veq}, it turns out that $\operatorname{span}(B)$ is also an invariant subspace of $F$, which means we can restrict the stability analysis of \eqref{Jacobianofmicronewcoord} to physically relevant perturbations by computing the matrix representation of $F: \operatorname{span}(B) \rightarrow \operatorname{span}(B) $ under a change of base to $B$ , denoted by $F^B_B$.

\begin{Lemma}
$\operatorname{span}(B)$ is an invariant subspace of $F$. Restricting $F$ to $\operatorname{span}(B)$ and changing base to $B$ gives
\begin{equation}
F_B^B = \begin{pmatrix}
0_{2N\times 2N} & \begin{matrix}
I_{2N-2 }\\
- 1^T_{N-1}\otimes I_2
\end{matrix} & 0_{2N\times 2} \\[6mm]
\lceil G(\hat{x})\rceil & - I_{N-1} \otimes 2\beta(m_0\otimes m_0^T)  & 0_{2N-2\times 2} \\[2mm]
0_{2\times 2N} & 0_{2 \times 2N-2} & -2\beta(m_0\otimes m_0^T)
\end{pmatrix}.
\label{FBB}
\end{equation}
\\
\proof
We compute the image of base vectors $b_i$ under $F$.
First, consider $b_i$ with a non-zero entry in a spatial coordinate ($i=1,\dots,2N$). Then
$Fb_i = (0_{2N}, g_{\cdot,i}, 0_2)^T$ consists of the $i$-th column of $G(\hat{x})=(g_{i,j})_{ij}$. Since the sum of each column of $G(\hat{x})$ in \eqref{hessian} is zero, we can rewrite this in terms of $B$ as
\begin{equation*}
 (0_{2N}, g_{\cdot,i}, 0_2)^T
=
\sum_{l=1}^{2N-2} g_{l,i} b_{2N+l}\,,
\end{equation*}
which gives the first block column of \eqref{FBB}.
Next, consider $b_i$ with a non-zero entry in a velocity component $k$ and assume $i$ to be odd $(i=2N+(2k-1)\,,\, k=1,\dots,N-1)$. Then we can compute $F b_i$ to get
\begin{equation*}
\begin{pmatrix} \delta^k_{2N-2} \\-1 \\0\\ \hdashline \\[-3mm]
-\begin{pmatrix}
  \frac{N-1}{N} \!& \!\dots \!&\!-\frac{1}{N} \\
  \vdots \!&\! \ddots \!&\! \vdots \\
 -\frac{1}{N} \!&\! \dots \!&\! \frac{N-1}{N}
  \end{pmatrix}
  \otimes 2\beta (m_0\otimes m_0^T)
   \begin{pmatrix}
 \delta^k_{2N-2} \\ -1\\0
  \end{pmatrix} \vspace{2mm}
 \\ \hdashline \\[-3mm]
 0_2
 \end{pmatrix}
 =
 \begin{pmatrix}
 \delta^k_{2N-2} \\ -1 \\ 0 \\ \hdashline
 0_{2(k-1)}\\
 -2\beta (m_0\otimes m_0^T)_{\cdot, 1}
 \\
 0_{2(N-k)-2} \\
 2\beta (m_0\otimes m_0^T)_{\cdot, 1}
 \vspace{2mm}
 \\
 \hdashline \\[-3mm] 0_2
 \end{pmatrix}
\end{equation*}
since, in the second block, the multiplication with $(\delta^k_{2N-2},-1,0)^T$ first gives $\frac{N-1}{N}+\frac{1}{N}$ in row $k$, $-\frac{1}{N}-\frac{N-1}{N}$ in row $2N-1$ and zero in all other rows, which is then multiplied with $-2\beta(m_0\otimes m_0^T)$ to select its first column twice (as $i$ is odd), indicated by the subscript $()_{\cdot,1}$.
This rewrites in terms of $B$ as
\begin{equation*}
Fb_i = b_{i-2N} - b_{2N-1}  - 2\beta m_{0,1}^2 b_{i} -2\beta m_{0,1}m_{0,2} b_{i+1}.
\end{equation*}
Analogously, one obtains for even $i$ $(i=2N+2k\,,\, k=1,\dots,N-1)$
\begin{equation*}
Fb_i = b_{i-2N} - b_{2N}  - 2\beta m_{0,1} m_{0,2} b_{i-1} -2\beta m_{0,2}^2 b_{i},
\end{equation*}
involving the second column of $-2\beta(m_0\otimes m_0^T)$.
The change of base for vectors of $B$ with non-zero entry in the mean velocity component is trivial.
Having shown that all images of base vectors can be expressed in terms of $B$, we showed the invariant subspace property and obtain the block structure of $F_B^B$ from our computations.
\QEDA
\label{cor-FBB}
\end{Lemma}

We now derive two corollaries on the linear structure of $F^B_B$, which will be essential to establish our stability result. To do so, we first rewrite \eqref{FBB} as
\begin{equation}
F_B^B =: \begin{pmatrix}
H & 0_{4N-2 \times 2} \\
0_{2 \times 4N-2} & -2\beta(m_0\otimes m_0^{T})
\label{def-H}
\end{pmatrix},
\end{equation}
since the two blocks separate. We first investigate the eigenspace associated to the zero eigenvalue of $F^B_B$. Some arguments are reminiscent of ideas in \cite[Lemma 3.2]{ABCV}.

\begin{Lemma} 
Suppose that $G(\hat{x})$ satisfies {\rm\bf (H1)-(H4)}, then $F_B^B$ does not posses a generalised eigenvector for eigenvalue zero.  \\
\proof Since the lower $2\times 2$ block of $F_B^B$ in \eqref{def-H} can be diagonalised with eigenvalues $0$ and $-2\alpha$, then we restrict our analysis to $H$. To simplify the notation, we will skip the dependency of the matrix $G(\hat{x})$ on the stationary state in this proof. If $z=(x,v)^T\in\eig(H,0)$ is a zero eigenvector, then $v = 0_{2N-2}$ and $x\in\eig(G,0)$ by \eqref{FBB}. Suppose $(u,w)^T$ is a generalized eigenvector, then
\begin{equation*}
H \begin{pmatrix}
u \\w
\end{pmatrix} = \begin{pmatrix}
x \\ 0_{2N-2}
\end{pmatrix}\,,\, x\in\eig(G,0).
\end{equation*}
which is equivalent to
\begin{subequations}
\begin{gather}
w = \lceil\mkern-2mu x\mkern-1mu\rceil,\label{gg1}\\
\sum_{i=1}^N x_{2i-1} = \sum_{i=1}^N x_{2i} = 0,  \label{gevproof-sumxzero}\\
\lceil\mkern-2mu Gu\mkern-1mu\rceil 
- I_{N-1}\otimes 2\beta(m_0\otimes m_0^T)  w = 0_{2N-2}. 
\label{gevproof-redeq}
\end{gather}
\end{subequations}
Due to the definition of $G$ in \eqref{hessian}, we get
\begin{equation}
\sum_{i=1}^N g_{2i-1,\cdot} = \sum_{i=1}^N g_{2i,\cdot} = 0_{2N}^T.
\label{zerosum}
\end{equation}
Adding two rows using \eqref{zerosum}, \eqref{gg1}, and \eqref{gevproof-sumxzero}, \eqref{gevproof-redeq} is equivalent to
\begin{equation}
Gu -  I_{N}\otimes 2\beta(m_0\otimes m_0^T) x = 0_{2N}. \label{gevproof-eq1}
\end{equation}
We left-multiply with $x$ and get
\begin{equation}
0=x^T G u = x^T (I_N \otimes 2\beta (m_0\otimes m_0^T )) x,
\end{equation}
since $G$ symmetric and $x\in\eig(G,0)$. This implies $x=[a_1 \dots a_N]^T\otimes m^{\perp} $ for some constants $a_i$ since the right-sided expression is quadratic and hence $x$ has to be particle-wise orthogonal to $m_0$ for the equation to hold ($m_0^\perp := (-m_{0,2}, m_{0,1})^T$).
On the other hand, \eqref{gevproof-sumxzero} implies that $x\in \operatorname{span}(w_3)$, since $w_3$, as given in \eqref{evdef-w3}, is the only zero eigenvector of $G$ whose component sums can be zero. Hence, all position vectors are tangential to a fixed mean velocity $m_0$, which is a contradiction to {\bf (H4)}.
\QEDA
\label{cor-nogev}
\end{Lemma}

Note that Lemma \ref{cor-nogev} is the main reason for this new change of variables \eqref{newcoord-veq}. Actually, the existence of generalized eigenvectors associated to the zero eigenvalue of the flock solution in the co-moving frame was shown in \cite[Remark 3.6]{ABCV}.
Next, we investigate the eigenvalues of $F^B_B$ under another hypothesis linking the eigenvectors of $G(\hat{x})$ with the orientationally arbitrary but fixed mean velocity $m_0$. Some arguments are reminiscent of ideas in \cite[Lemma 3.1]{ABCV}. 

\begin{Lemma} Suppose $G(\hat{x})$ satisfies {\rm\bf (H1)-(H4)} and together with $m_0$ satisfies {\rm\bf (H5)}, then $\dim(\eig(F_B^B,0))=4$. Furthermore, all non-zero eigenvalues $\mu_i$ of $F_B^B$ have negative real-part, $\mu_i\neq 0 \Rightarrow \re(\mu_i)<0$. \\
\proof
Lemma {\rm\ref{cor-nogev}} shows that there are no generalized eigenvectors for the zero eigenvalue of $F_B^B$. We will again skip the dependency of the matrix $G$ on the stationary state. Let us first show that $\dim(\eig(F_B^B,0))\geq 4$. One directly checks that  the basis elements $\{w_1,w_2,w_3\}$ of $\eig(G,0)$ define linearly independent eigenvectors in $Eig(F_B^B,0)$ by $z_i := (w_i , 0_{2N-2} , 0_{2} )^T$, $i=1,2,3$. As remarked in Lemma \ref{cor-nogev}, the lower $2\times 2$ block of $F^B_B$ can be diagonalised with eigenvalues $0$ and $-2\alpha$ and associated eigenvectors $z_4:=(0_{2N}, 0_{2N-2}, m_0^{\perp})^T$ and $(0_{2N},0_{2N-2}, m_0)^T$ respectively. 

Next, we prove the last statement and $\dim(\eig(F_B^B,0))\leq 4$ simultaneously.
First, note that any pair of eigenvalue $\mu$ and eigenvector $z=(x,v)^T$  satisfies $Hz = \mu z$, which gives
\begin{subequations}
\begin{gather}
v = \mu \lceil\mkern-2mu x\mkern-1mu\rceil \label{spectrumproof-vismux} \\
\sum_{i=1}^N x_{2i-1} = \sum_{i=1}^N x_{2i} = 0  \label{spectrumproof-sumxzero} \\
\lceil\mkern-2mu Gx\mkern-1mu\rceil - ( I_{N-1} \otimes 2\beta (m_0\otimes m_0^T))v = \mu v.  \label{spectrumproof-redeq}
\end{gather}
\end{subequations}
Equation \eqref{spectrumproof-redeq} rewrites as
\begin{equation}
\mu^2 \lceil\mkern-1mu x\mkern-1mu\rceil +  
(I_{N-1}\otimes 2\beta\mu (m_0\otimes m_0^T))
\lceil \mkern-1mu x\mkern-1mu\rceil = 
\lceil\mkern-2mu Gx\mkern-1mu\rceil.  \label{spectrumproof-redeq2}
\end{equation}
By the same argument used in Lemma \ref{cor-nogev}, \eqref{spectrumproof-redeq2} is equivalent to
\begin{equation*}
\mu^2 x +  ( I_{N} \otimes 2\beta\mu (m_0\otimes m_0^T))x = Gx,  \label{spectrumproof-nonredeq}
\end{equation*}
due to \eqref{spectrumproof-sumxzero} and the properties of $G$.
We left-multiply with $\bar{x}$ to obtain
\begin{equation*}
\mu^2 + 2\beta\mu \sum_{i=1}^N \scal{m_0}{x_{\substack{2i\text{-}1\\ \!\!\!\!\!2i}}}\overline{\scal{m_0}{x_{\substack{2i\text{-}1\\ \!\!\!\!\!2i}}}} = \bar{x}^T Gx.
\end{equation*}
assuming \textwlog the normalisation $\scal{x}{x}=\bar{x}^Tx=1$ and $\scal{\cdot}{\cdot}$ denoting the complex scalar product. Therefore, we deduce
\begin{equation*}
\mu=-A\pm\sqrt{A^2-B}
\text{ with }
A = -\beta \sum_{i=1}^N \left|\scal{m_0}{x_{\substack{2i\text{-}1\\ \!\!\!\!\!2i}}}\right|^2 
\,\text{and }\,
B = -\bar{x}^T Gx.
\end{equation*}
Note that $A\geq 0, B\geq 0$ since $\bar{x}^T Gx\leq 0$ for any $x$ due to {\rm\bf (H3)}.
We study three cases:\\
\underline{Case 1:} Suppose $A>0, B>0$. If $A^2-B\geq 0 \Rightarrow \mu<0$ by monotonicity. If $A^2-B<0$, its square root is purely imaginary and thus $\re(\mu_i)<0$.   \\
\underline{Case 2:} Suppose $A=0,B>0$. Then $\mu=\pm\sqrt{\bar{x}^TGx}\in i\rr$.
Since $A$ is a sum over squared moduli of $N$ complex numbers, $A = 0 \Leftrightarrow \scal{m_0}{x_{\substack{2i\text{-}1\\ \!\!\!\!\!2i}}}= 0$ for all $i$. Since $m$ is fixed, this holds if and only if $x$ is particle-wise orthogonal to $m_0$.
Now, inserting $A=0$ back into \eqref{spectrumproof-nonredeq} gives $Gx=\mu^2 x, \mu^2<0$ which implies that such an $x$ is also an eigenvector for $G$. This is a contradiction to {\rm\bf (H5)} and no such case exists.
\\
\underline{Case 3:} Suppose $A\geq0, B=0$. Then $x\in \eig(G,0)$ and either $\mu=-2A<0$ 
or $\mu=0$. For the second case, $v=0_{2N-2}$ follows from \eqref{spectrumproof-vismux}, which implies $z\in\operatorname{span}(z_1,z_2,z_3)\subset \eig(F_B^B,0)$ by {\rm \bf (H2)}.

Together, this shows that
 \begin{subequations}
 \begin{gather}
 \mu_i \neq 0 \Rightarrow \re(\mu_i) < 0. \\
 A\geq 0 , B=0,  \mu=0 \Rightarrow z\in\operatorname{span}(z_1,z_2,z_3)\subset\eig(F_B^B,0). \label{secondconcl}
 \end{gather}
\end{subequations}
To complete the study of $\eig(F_B^B,0)$ suppose now that $\mu=0$. Then $A^2-B\geq 0$ and $0=-A+\sqrt{A^2-B}\Leftrightarrow A^2 = A^2-B \Leftrightarrow B=0$, which takes us back \eqref{secondconcl}. Therefore, no further eigenvector to eigenvalue zero exists, which completes the proof.
\QEDA
\label{cor-spec}
\end{Lemma}

Let us remark that the 4 basis eigenvectors associated to the zero eigenvalue of $F_B^B$ represent linearized flow within the set of stationary flock solutions $\tilde{\mathcal{Z}}_F$ in \eqref{flocks2}.
From the change of variables \eqref{newcoord} and Lemma \ref{cor-newcoord}, we know that the family of flock solutions is written as $\hat{Q}=(\hat{x}, 0_{2N-2}, m_0)^T$
with $ |m_0|^2=\tfrac{\alpha}{\beta}$.  $G(\hat{x})$ is invariant under translation and spatial rotation, and thus $z_1$, $z_2$, and $z_3$ are the linearized flow along these invariances. The fourth eigenvector $z_4$ represents linearized rotation of the mean velocity $m_0$, compared to \cite[Remark 3.6]{ABCV}.

%%%%%%%%%%%%%%%%%%%%%%%%%%%%%%%%%%%%%%%%%

\section{Proof of the main result}

Changing coordinates to \eqref{newcoord} eliminates the linear translation at constant speed and flock solutions become constant solutions which we identify with their profile and velocity at $t=0$. The image of $\mathcal{Z}_F$ under \eqref{newcoord} is $\tilde{\mathcal{Z}}_F$ in \eqref{flocks2}. $\tilde{\mathcal{Z}}_F$ is a manifold of stationary states of the flow of system \eqref{newcoord-veq}. We aim to show that $\tilde{\mathcal{Z}}_F$ is a normally hyperbolic invariant manifold according to \cite{HPS} for the time flow map of \eqref{newcoord-veq} restricted on $\operatorname{span}(B)$, which we denote by $\mathbb{F}^t_B$.
First, $\tilde{\mathcal{Z}}_F$ is trivially invariant, $\mathbb{F}^t_B(\tilde{\mathcal{Z}}_F)=\tilde{\mathcal{Z}}_F$ for all $t\in\rr$, as it is entirely stationary by Lemma \ref{cor-newcoord}.
Choose $Q^*\in \tilde{\mathcal{Z}}_F$ such that {\bf (H5)} is fulfilled. Since eigenvectors are distinct, {\bf (H5)} holds true for a neighbourhood of $Q^*$. There is no tangential flow along $\tilde{\mathcal{Z}}_F$, and in particular
\begin{equation}\label{flowddt}
\frac{\dd}{\dd t} D_{Q}\mathbb{F}^t_B (Q^*) = F^B_B \left[ D_Q\mathbb{F}^t_B(Q^*)\right].
\end{equation}
Corollaries \ref{cor-nogev} and \ref{cor-spec} tell us, that $F_B^B$ decomposes into the direct sum of $\eig(F_B^B,0)$, which equals the tangential space of $\tilde{\mathcal{Z}}_F$ at $Q^*$, and a subspace of eigenvectors whose eigenvalues have negative real part and span the non-tangential subspace of $\tilde{\mathcal{Z}}_F$ at $Q^*$. It is important to stress, that this splitting is preserved under the flow $\mathbb{F}^t_B$, since $Q^*$ is stationary.
Adapting from \cite{HPS}, all prerequisites for $\tilde{\mathcal{Z}}_F$ to be a normally hyperbolic invariant manifold are given. In the language of dynamical systems theory, the tangent bundle of $\rr^{4N}$ at $Q^*$ has an $\mathbb{F}_B$-invariant splitting into tangential and stable sub-bundles.
A necessary inequality between the tangential flow and the stable decay rate is trivially given in our case, since flow along $\tilde{\mathcal{Z}}_F$ is zero and as a consequence, we have a spectral gap.
Note, that though $\tilde{\mathcal{Z}}_F$ is not compact we can easily select a compact sub-manifold around $Q^*$ or perform a symmetry reduction and eliminate the translational invariance by further reducing $\operatorname{span}(B)$.

We can hence apply the Theorem of Hadamard and Perron as in \cite[Theorem 4.1]{HPS}, which, adapted to our situation, gives the following:
There exists locally an $\mathbb{F}_B$-invariant stable smooth manifold in a neighbourhood of $\tilde{\mathcal{Z}}_F$, where the splitting into stable and stationary of the linearized flow is preserved.
The space between both manifolds is laminated by invariant fibers tangent to the stable subspace at $Q^*$, which are curves of the nonlinear flow towards $\tilde{\mathcal{Z}}_F$.
Hence, for a small perturbation of $Q^*$ within that laminated space, we are sure to be sitting on some fiber, which might have $Q^*$ as its base point on $\tilde{\mathcal{Z}}_F$ or a point in a neighbourhood of $Q^*$. In either case, there is transport along the fiber towards $\tilde{\mathcal{Z}}_F$, for which an exponential decay estimate holds.
This estimate is obtained by closing $\eqref{flowddt}$ on the stable sub-bundle.
In other words, there exists locally a decomposition into fast and slow (stationary, in our case) coordinates of the nonlinear flow, such that the dynamics is characterised by the linearized flow.
Going back to original coordinates \eqref{micro}, we obtain the statements in Theorem \ref{theo-main} and Corollary \ref{cor-main}.

\begin{remark}
The precise meaning of "relaxation in time" towards a different flock solution in the statement of Theorem \ref{theo-main} is given by the notion of stable hyperbolicity of the family of flock solutions in coordinates \eqref{newcoord}, as used in the proof. In the original coordinates, there is no notion of ''closeness'' between elements of \eqref{flocks},
since e.g. a change of orientation will cause two flocks to arbitrarily diverge in euclidean norm as $t\rightarrow \pm \infty$.
Relaxation thus means, that the trajectory of a perturbed flock will approach a trajectory of another flock solution with typically different orientation and similar spatial configuration.
Theorem \ref{theo-main} also tells us that the manifold is stable towards small perturbations of the flow $\mathbb{F}_B$, though that case is not considered here. Note, that this allows for flock stability under small perturbations of the potential.
\end{remark}

\section{Numerical experiments}
In this section, we numerically investigate the hypothesis of Theorem \ref{theo-main} for a number of commonly used interaction potentials and study the nonlinear stability of flocks in terms of their polarization. 

The key advantage of Theorem \ref{theo-main} is that essentially linear (asymptotic except symmetries) stability of the spatial configuration of a flock in the first-order aggregation model \eqref{aggregation} is inherited nonlinearly by the family of flock solutions to the second-order model \eqref{micro}.   
In order to demonstrate that the stability theorem applies to a given potential, {\bf (H1)-(H4)} have to hold for $\hat{x}$ being a stationary state of \eqref{aggregation}. 
We compute a numerical approximation $\{\hat{x}\}$ of $\hat{x}$ by a standard fourth-order Runge-Kutta time integration of \eqref{aggregation} with random initial data and an  interaction potential scaled with a factor $D$ to relate the amplitude of the right-hand side of \eqref{aggregation} with the accuracy of the solver:
\begin{equation*}
W_D(r) := D\, W(r). 
\label{potscaling}
\end{equation*}
Note, that $\hat{x}$ is independent of $D$. We hence obtain a spatial configuration for which {\bf (H1)} holds numerically. Let us mention that {\bf (H4)} is clearly satisfied for all potentials considered here. 
Hypotheses {\bf (H2)-(H3)} require the Jacobian $G(\{\hat{x}\})$ to posses negative eigenvalues except for three zero eigenvalues associated to translation and rotation.
As $G(\{\hat{x}\})\in \rr^{2N\times 2N}$ is symmetric it is sufficient to apply standard methods to compute its spectrum, at least for the number of particles we consider.
Let $(\mu_1,\mu_2;\mu_3)$ denote the eigenvalues of $G(\{\hat{x}\})$ associated to translation and rotation and $\mu_i, i\geq 4$ all other eigenvalues sorted in descending order.  
In Table \ref{tab-ev} we show $\mu_4$ for different potentials and numbers of particles, together with $|\mu_3|$ to illustrate the numerical accuracy of $\{\hat{x}\}$ and the scaling parameter $D$. For the comparability of the potentials, $D$ is set such that $\mu_{2N}=1$ for all cases shown. 
The potentials used in Table \ref{tab-ev} are:
\begin{equation*}
\begin{tabular}{l | l|}
& \multicolumn{1}{l}{$W(r) = V(r) - CV(r/\ell) $} \\
\cline{1-2}
\multicolumn{1}{l|}{Morse}  & $V(r) = -e^{-r}$\\
\multicolumn{1}{l|}{Quasi-Morse }  & $V(r) = -\frac{1}{2\pi}K_{0}(kr)$\\
\multicolumn{1}{l|}{Generalized Morse}  & $V(r) = -e^{-\frac{r^p}{p}}$\\
\hline \hline
\multicolumn{1}{l|}{Log-Newtonian} & \multicolumn{1}{ l }{$W(r) = r^2 - \log(r)$} \\
\end{tabular}
\end{equation*}
where $K_0$ is the modified Bessel function of second kind (for detailed discussions of the respective potentials see \cite{CHM, Carrillo2013, CDMBC, FHK}).

\begin{table}
\def\firstcolumnlength{.35\textwidth}
\centering
\begin{tabular}{cc|c||c|c|}
\cline{3-5}
& & \multicolumn{3}{ c| }{Eigenvalue structure of $G$} \\ \cline{2-5}
& \multicolumn{1}{|c|}{N} & $\mu_4$ & $|\mu_3|$ & $D$  \\ \cline{1-5}
\multicolumn{1}{ |c| }{\multirow{4}{*}{\begin{minipage}{\firstcolumnlength}\centering Morse potential \\ $C\!=\!\frac{10}{9}, \ell\!=\!\frac{3}{4}$ \end{minipage}} } &
\multicolumn{1}{ |c| }{25} & -6.9579e-5 &  5.3279e-9 &  9.5290     \\ \cline{2-5}
\multicolumn{1}{ |c  }{}                        &
\multicolumn{1}{ |c| }{40} & -3.9449e-4 & 3.0435e-8 &  9.7423     \\ \cline{2-5}
\multicolumn{1}{ |c  }{}                        &  
\multicolumn{1}{ |c| }{70} & -2.3600e-4 & 6.8493e-9 & 10.0008   \\ \cline{2-5}
\multicolumn{1}{ |c  }{}                        &  
\multicolumn{1}{ |c| }{100} & -2.7069e-5 & 2.7073e-8 & 10.0720 \\ \cline{1-5}
\multicolumn{1}{ |c| }{\multirow{4}{*}{\begin{minipage}{\firstcolumnlength}\centering Quasi-Morse potential \\ $C\!=\!\frac{10}{9}, \ell\!=\!\frac{3}{4}, k\!=\!\frac{1}{2}$ \end{minipage}} } &
\multicolumn{1}{ |c| }{25} & -6.7374e-4 & 1.4273e-11 & 39.8647     \\ \cline{2-5}
\multicolumn{1}{ |c  }{}                        &
\multicolumn{1}{ |c| }{40} & -1.6636e-4 & 3.5242e-12 & 39.0230      \\ \cline{2-5}
\multicolumn{1}{ |c  }{}                        &  
\multicolumn{1}{ |c| }{70} &  -0.0024 & 9.1205e-12 & 38.5022   \\ \cline{2-5}
\multicolumn{1}{ |c  }{}                        &  
\multicolumn{1}{ |c| }{100} & -0.0023 & 4.0528e-11 & 38.3616 \\ \cline{1-5}
\multicolumn{1}{ |c| }{\multirow{4}{*}{\begin{minipage}{\firstcolumnlength}\centering Generalized Morse potential \\ $C\!=\!\frac{10}{9}, \ell\!=\!\frac{3}{4}, p=1.25$ \end{minipage}} } &
\multicolumn{1}{ |c| }{25} &-5.0193e-5 &  4.5907e-8 & 7.6429      \\ \cline{2-5}
\multicolumn{1}{ |c  }{}                        &
\multicolumn{1}{ |c| }{40} & -2.7459e-4 & 1.8389e-8 &  7.8653     \\ \cline{2-5}
\multicolumn{1}{ |c  }{}                        &  
\multicolumn{1}{ |c| }{70} & -1.2202e-4 & 1.7075e-10 & 7.9694   \\ \cline{2-5}
\multicolumn{1}{ |c  }{}                        &  
\multicolumn{1}{ |c| }{100} & -8.5882e-6 & 1.4491e-10 & 7.9813 \\ \cline{1-5}
\multicolumn{1}{ |c| }{\multirow{4}{*}{\begin{minipage}{\firstcolumnlength}\centering Log-Newtonian potential 
 \end{minipage}} } &
\multicolumn{1}{ |c| }{25} & -3.6270e-4 &  4.0557e-6 & 0.25     \\ \cline{2-5}
\multicolumn{1}{ |c  }{}                        & 
\multicolumn{1}{ |c| }{40} &-1.8307e-8  & 8.0974e-9 & 0.25      \\ \cline{2-5}
\multicolumn{1}{ |c  }{}                        &  
\multicolumn{1}{ |c| }{70} & -0.0022 & 6.1172e-7 & 0.25   \\ \cline{2-5}
\multicolumn{1}{ |c  }{}                        &  
\multicolumn{1}{ |c| }{100} & -0.0029 & 7.1695e-6 & 0.25 \\ \cline{1-5}
\end{tabular}
\caption{Numerical results for the largest non-zero eigenvalue $\mu_4$ of the Jacobian $G(\{\hat{x}\})$ for varying potentials $W$ and number of particles $N$. The numerical value of the zero eigenvalue $\mu_3$ is shown as a measure of the exactness of $\{\hat{x}\}$. The free scaling parameter $D$ is set to normalize the smallest eigenvalue to $\min(\mu_i)=-1$.
As $\mu_4$ is negative, hypothesis {\bf (H2)-(H3)} of Theorem \ref{theo-main} are numerically validated for the cases shown.  
}
\label{tab-ev}
\end{table}
 
It can be seen that for all cases considered, eigenvalues other than $\mu_1,\mu_2,\mu_3$ are negative and hence Theorem \ref{theo-main} applies.
The difference between the numerical eigenvalue $\mu_3$, which is equal to zero in theory, and $\mu_4$ is of order four or higher, which shows that the negativity of $\mu_4$ is not a numerical artefact. However, we also observe that smaller gaps occasionally occur for some configurations ($N=40$, Log-Newtonian) due to additional quasi-symmetries in $\{\hat{x}\}$. 
We remark, that the eigenvalue computations of Table \ref{tab-ev} require a highly accurate computation of $\{\hat{x}\}$, which is here obtained with a scaling of $D=500$ and a time interval $T=[0,500]$.   

\begin{figure}[h!]
\subfloat[Statistical comparison of initial polarization \newline vs.\@ minimal polarization in time evolution.]{\includegraphics[keepaspectratio=true,
width=.5\textwidth]{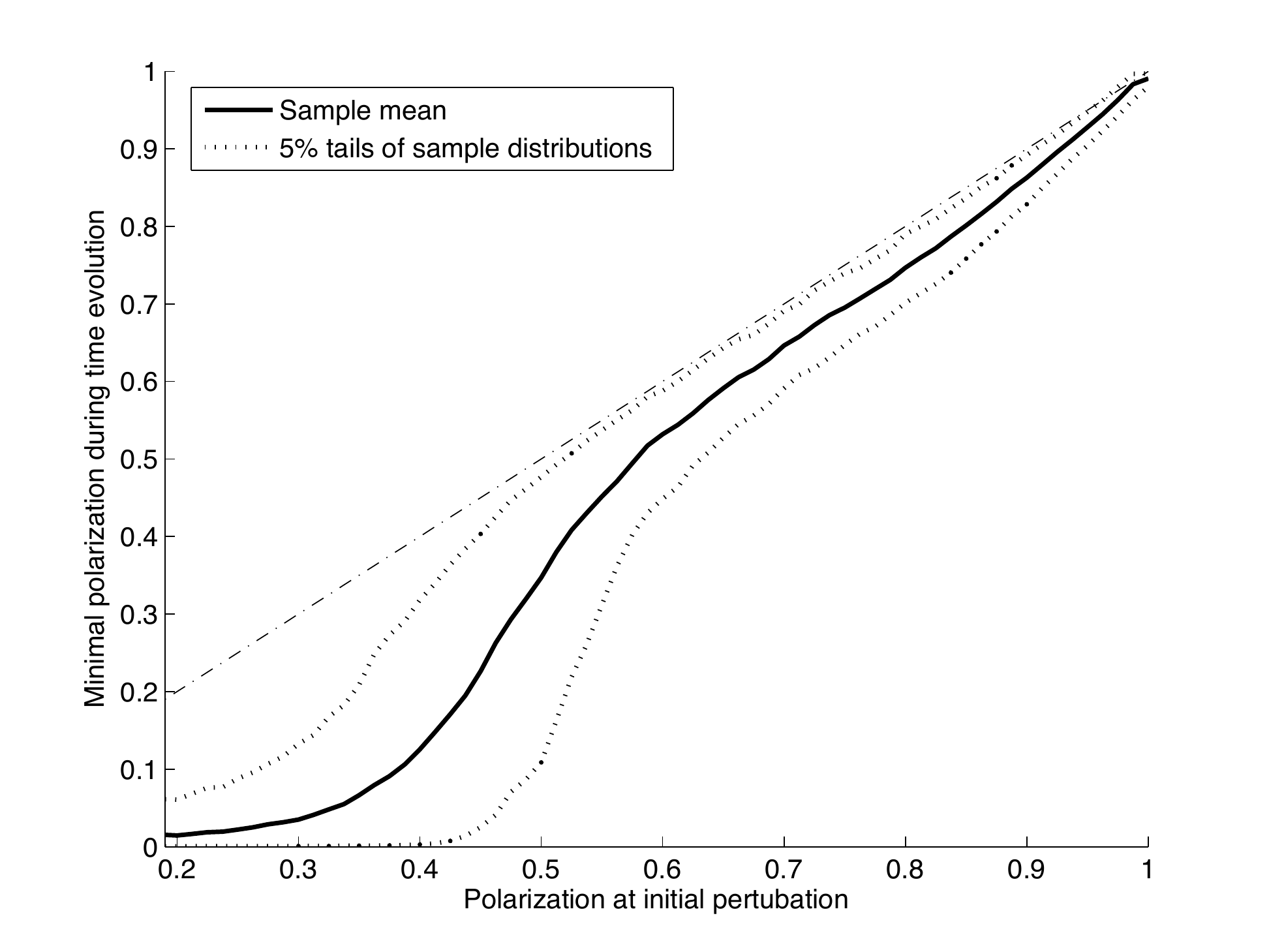}}
\subfloat[Statistics of the mean-$\ell_2$ perturbation per particle in $(x,v)$ for the initial perturbations of (a).]{\includegraphics[keepaspectratio=true,
width=.5\textwidth]{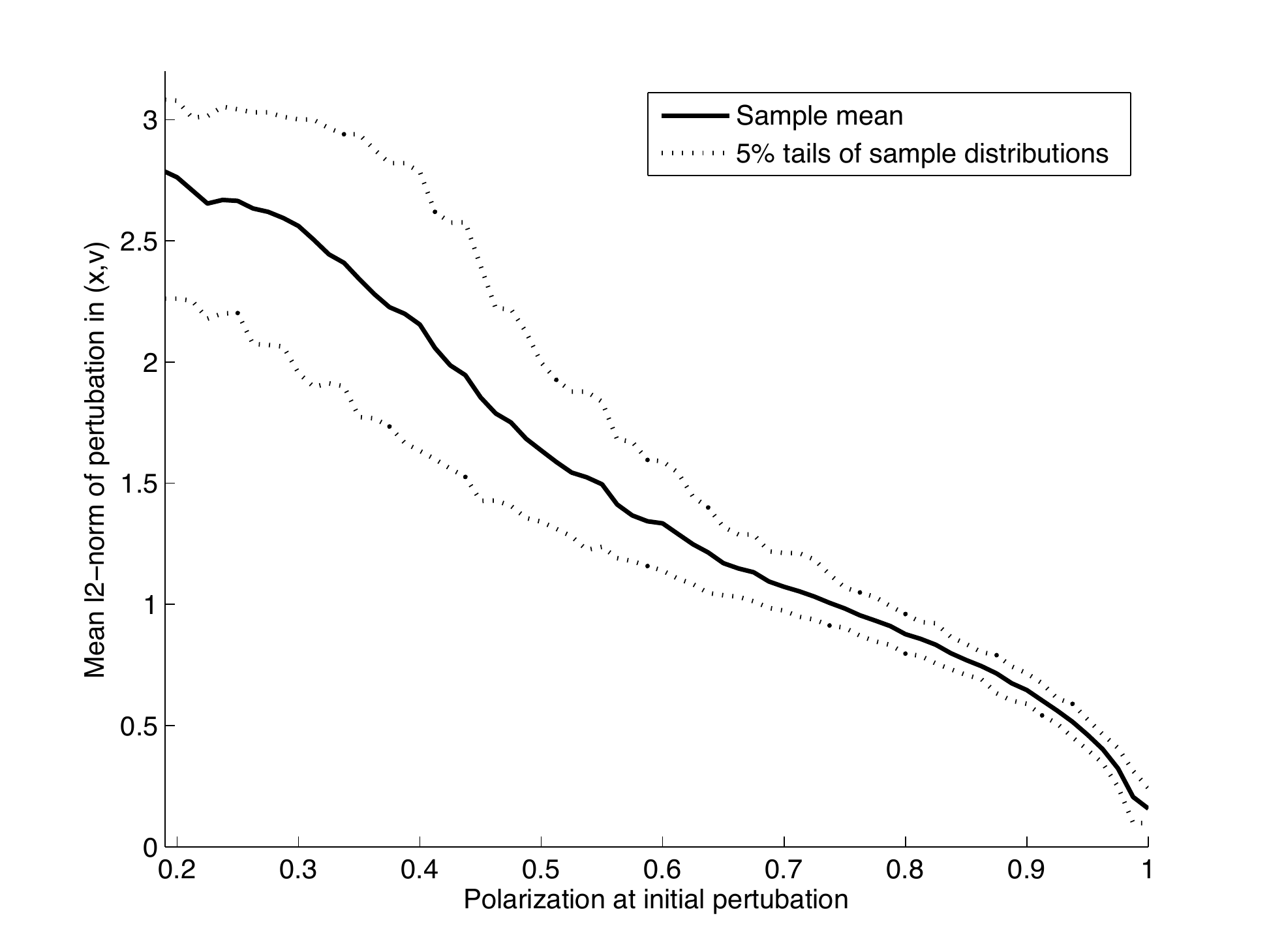}}
\caption{Numerical findings of a sequence of simulations of system \eqref{micro} with perturbed initial data. (a): Histogram of initial polarization of the perturbed flock against the corresponding minimal polarization observed during the subsequent time evolution. (b): Supporting statistics of the initial polarization against the average $\ell_2$-norm of perturbation per particle that induce the respective polarisation change. From (a) ($-$: arithmetic mean, $-\,-$: 5\% tails of statistical distribution, Morse potential as in Table \ref{tab-ev}, $\alpha=1,\beta=5, \lambda=50$, $N=100$ particles, 25000 simulations}
\label{fig-pol}
\end{figure}

Next, we aim to study the stability of a flock in the nonlinear dynamics of \eqref{micro} for the Morse potential under uniformly distributed perturbations in space and velocity. Theorem \ref{theo-main} provides a theoretical bound on the decay rate into equilibrium, which is very small as seen in Table \ref{tab-ev}.  
For random perturbations it is natural to assume that the flock is ''more stable'' than the theorem asserts.
To study whether a flock remains a flock after an initial perturbation and the subsequent dynamics of \eqref{micro}, we use the notion of polarization of a group of individuals, which is defined as
\begin{equation*}
\operatorname{Pol}(\{(x_i,v_i)\}_{i=1}^N) := \frac{1}{N}\sum_i \frac{v_i}{\norm{v_i}},
\end{equation*}
for $N$ agents with phase-space coordinates at $(x_i,v_i)\in\rr^4$, see e.g. \cite{couzinKJRF} for instance. 
For an aligned flock, the polarization equals one whereas e.g. a rotating mill possesses zero polarization. 
We perform a sequence of simulations of \eqref{micro} with randomly perturbed initial data. 
The stationary flock solution $z=\{(\hat{x}_i, m_0)\}_{i=1}^N), \norm{m_0}=\frac{\alpha}{\beta}$ 
is perturbed by $\{(\Delta x_i, \Delta v_i)\}_{i=1}^N$ with $(\Delta x_i, \Delta v_i)$ sampled from the uniform distribution on $[-\frac{a}{2},\frac{a}{2}]\times [-\frac{a}{2},\frac{a}{2}]$ and varying perturbation strength $a$.
Stability of the flock in the time evolution $[0,T]$ can estimated as follows: If the support of the solution remains bounded and the polarization is equal or larger than the change in polarization induced by the initial perturbation, is it reasonable to assume that the ensemble will again tend to form a flock solution. If however, the polarization drops to significantly low or near zero values, the particles have evolved into a chaotic state whose convergence behaviour as $T\rightarrow \infty$ is unknown. 
In Figure \ref{fig-pol}(a), we plot the initial polarization induced by the perturbation against the minimal perturbation that occurred during a time integration of \eqref{micro} with $T=100$. The figure is based on 25000 simulations with uniformly distributed perturbations strength $a$ and shows the statistical mean and the $5\%$ tails of the statistical distribution.  
In Figure \ref{fig-pol}(b), we show the correlation between the polarization induced by the initial perturbation and the averaged $l_2$-norm of the perturbation per particle $\frac{1}{N}\sum_i \norm{(\Delta x_i,\Delta v_i)}_{2}$, using the same statistical measures. 
The potential parameters are identical to Table \ref{tab-ev} and the number of particles is $N=100$.   
We observe from Figure \ref{fig-pol}, that the probability of loosing additional polarisation in the flock nonlinearly increases with the strength of the initial perturbation. Below a critical threshold of $\approx 0.6$, the flock will most likely drop polarization with an increasing probability of near zero values (e.g. chaotic state).    
From Figures \ref{fig-pol}(b) we observe that a perturbation-induced polarization of $0.6$ is on average generated by a strong average perturbation in $(x,v)$ of $\approx 1.34$ per particle 
hence the nonlinear stability observed in this experiment well exceeds the validity regime of Theorem \ref{theo-main}. 

\section{Conclusions}
In this paper, we have shown that flock solutions are fully nonlinearly stable under small enough position-velocity perturbations, in the sense that a perturbed flock configuration will relax to another flock configuration in $\tilde{\mathcal{Z}}_F$. Stability is inherited from the first-order particle system to the second-order model.
Hence the flock might slightly rotate its position, translate its centre of mass, change the orientation of its mean velocity, or any combination of the three, but it will remain a flock solution.
This result has been established for general interaction potentials $W$ satisfying assumptions ${\bf (H1)-(H4)}$. 
Our findings are completely consistent with one's expectations from numerical experiments and, to the best of our knowledge, is the first of its kind for general flocks in second-order models.
The task of theoretically proving nonlinear stability results of steady solutions for the first-order aggregation equation remains largely open, and is inevitably complex given the extensive variety of stationary states (see e.g. \cite{BCLR2}). We therefore numerically investigated the spectrum of the Jacobian $G(\hat{x})$ evaluated at a stationary state $\hat{x}$ for some commonly used potentials, and found them to posses the desired eigenvalue structure. Flocks for Morse-type potentials are linearly unstable in the co-moving frame second-order model (\cite{ABCV}), but nonlinearly asymptotically stable in the sense of our main result in Theorem \ref{theo-main}.
Finally, we investigated the nonlinear stability of flocks under uniformly distributed perturbations, and numerically validated the nonlinear stability of flocks in the sense of polarization.

\subsection*{Acknowledgements}
JAC acknowledges support from projects MTM2011-27739-C04-02,
2009-SGR-345 from Ag\`encia de Gesti\'o d'Ajuts Universitaris i de Recerca-Generalitat de Catalunya, and the Royal Society through a Wolfson Research Merit Award. JAC, YH, and SM were supported by Engineering and Physical Sciences Research Council (UK) grant number EP/K008404/1.

\end{document}